\documentclass[12pt]{article}

\newtheorem{theorem}{Theorem}

\newtheorem{definition}{Definition}

\newcommand{\JK}[2]{
\begin{eqnarray*}
J & = & \{ #1 \} \\
K & = & \{ #2 \}
\end{eqnarray*}
}

\usepackage{epsfig}
\usepackage{url}
\usepackage{amssymb}
\usepackage{amsmath}
\usepackage{amsfonts}

\def\Dbar{\leavevmode\lower.6ex\hbox to 0pt{\hskip-.23ex \accent"16\hss}D}

\def\bZ{{\mbox{\bf Z}}}

\def\naf{{\mbox{\rm NAF}}}
\def\paf{{\mbox{\rm PAF}}}
\def\psd{{\mbox{\rm PSD}}}

\begin{document}

{\bf\LARGE
\begin{center}
Some new periodic Golay pairs
\end{center}
}

{\Large
\begin{center}
Dragomir {\v{Z}}. {\Dbar}okovi{\'c}
\footnote{University of Waterloo, Department of Pure Mathematics
and Institute for Quantum Computing,
Waterloo, Ontario, N2L 3G1, Canada,
e-mail: \url{djokovic@uwaterloo.ca}},
Ilias S. Kotsireas
\footnote{Wilfrid Laurier University,
Department of Physics \& Computer Science,
Waterloo, Ontario, N2L 3C5, Canada,
e-mail: \url{ikotsire@wlu.ca}}
\end{center}
}
\begin{abstract}
Periodic Golay pairs are a generalization of ordinary Golay pairs.  They can be used to construct Hadamard matrices.
A positive integer $v$ is a (periodic) Golay number if there exists a (periodic) Golay pair of length $v$. Taking into the account the results obtained in this note and yet unpublished new result \cite{DK::2014}, there are only seven known periodic Golay numbers which are definitely not Golay numbers, namely $34,50,58,68,72,74,82$. We construct here periodic Golay pairs
of lengths $74,122,164,202,226$. It is apparently unknown whether $122,164,202,226$ are Golay numbers.
The smallest length for which the existence of periodic Golay pairs is undecided is now $90$.
\end{abstract}

\section{Preliminaries} \label{prelim:sec}

Let $v$ be a positive integer and $\bZ_v=\{0,1,\ldots,v-1\}$
the ring of integers modulo $v$. Let $A=a_0,a_1,\ldots,a_{v-1}$ be a binary, i.e., $\{\pm1\}$-sequence of length $v$. The periodic autocorrelation function $\paf_A:\bZ_v\to\bZ$ and the nonperiodic autocorrelation function $\naf_A:\bZ\to\bZ$ of $A$ are defined by
\begin{eqnarray*}
\paf_A(s) &=& \sum_{i=0}^{v-1} a_i a_{i+s\pmod{v}}, \\
\naf_A(s) &=& \sum_{i\in\bZ} a_i a_{i+s}.
\end{eqnarray*}
(In the nonperiodic case we set $a_i=0$ when $i<0$ or $i\ge v$.)
Note that
\begin{equation} \label{paf-naf}
\paf_A(s)=\naf_A(s)+\naf_A(v-s), \quad s\in\bZ_v.
\end{equation}
A {\em Golay pair} is an ordered pair $(A,B)$ of binary sequences of length $v$ such that $\naf_A(s)+\naf_B(s)=0$ for
$s\ne0$.
Similarly, a {\em periodic Golay pair} is an ordered pair
$(A,B)$ of binary sequences of length $v$ such that
$\paf_A(s)+\paf_B(s)=0$ for $s\ne0$.
It follows from (\ref{paf-naf}) that any Golay pair is also
a periodic Golay pair.
Periodic Golay pairs are also known as {\em complementary
binary sequences}. They can be used to construct Hadamard matrices of order $2v$ (see subsection \ref{HadMat} below). 
Thus, if $v>1$ is a periodic Golay number then $v$ must be even.

Periodic Golay pairs can be viewed as a particular case
of supplementary difference sets (SDS). Let us recall some facts
which will be used in our constructions, see e.g.
\cite[Proposition 1]{DK:compression:2013}.

First we have to recall the definition of SDS.
Let $k_1,\ldots,k_t$ be positive integers and $\lambda$ an integer such that $\lambda(v-1)=\sum_{i=1}^t k_i(k_i-1)$,
and let $X_i$ be a subset of $\bZ_v$ with cardinality
$|X_i|=k_i$, $i\in\{1,2,\ldots,t\}$.

\begin{definition} \label{def-sds}
We say that $(X_1,\ldots,X_t)$ are {\em supplementary difference
sets} with parameters $(v;k_1,\ldots,k_t;\lambda)$, if for every nonzero element $c\in\bZ_v$ there are exactly $\lambda$
ordered pairs $(a,b)$ such that $a-b=c \pmod{v}$ and
$\{a,b\}\subseteq X_i$ for some $i\in\{1,2,\ldots,t\}$.
\end{definition}

For convenience, we shall refer to the sets $X_1,\ldots,X_t$
as the {\em base blocks} of this SDS. To any SDS with parameters $(v;k_1,\ldots,k_t;\lambda)$ we attach an additional parameter
$n$ defined by $n = k_1 + \cdots + k_t - \lambda$.
(SDSs with $t=1$ are known as cyclic difference sets.)

There is a bijection from the set of all binary sequences of
length $v$ to the set of all subsets of $\bZ_v$ which assigns
to the sequence $A$ the subset $\{i\in\bZ_v:a_i=-1\}$. If
$(A,B)$ is a periodic Golay pair, then the corresponding
pair of subsets $(X,Y)$ of $\bZ_v$ is an SDS whose parameters
$(v;r,s;\lambda)$ satisfy the equation $v=2n$.
(Recall that $n=r+s-\lambda$ in this case.)
The converse is also true.
Moreover, if $a=v-2r$ and $b=v-2s$ then $a^2+b^2=2v$. In particular, $v$ must be a sum of two squares.

\subsection{Periodic Golay pairs and Hadamard matrices}
\label{HadMat}

One of the reasons that periodic Golay pairs are useful is that they can be used to
construct Hadamard matrices. Let us detail the relevant construction. Suppose that two binary sequences $A, B$ of length $v$
that form a periodic Golay pair are given.
Then one can construct the associated $v \times v$
circulant matrices $C_A, C_B$ whose first rows are
the sequences $A, B$ respectively.
These two circulant matrices satisfy the matrix equation
$$
    C_A C_A^t + C_B C_B^t = (2v) I_v
$$
where $t$ denotes transposition and $I_v$ denotes the $v \times v$ unit matrix.
Using these two circulant matrices, a Hadamard matrix of order $2v$ can be constructed as
\begin{equation} \label{Konstrukcija}
    H_{2v} =
    \left[
    \begin{array}{c|c}
    C_A & C_B \\
    \hline
    -C_B^t & C_A^t \\
    \end{array}
    \right] .
\end{equation}
We refer the reader to \cite{SY:1992} for more details
and additional Hadamard matrix constructions.

\section{Known Golay and periodic Golay numbers}
\label{Known}

If $\alpha,\beta,\gamma$ are nonnegative integers, then it
is well known that $2^\alpha10^\beta26^\gamma$ is a Golay number.
No additional Golay numbers have been found so far.
Moreover, there are no other Golay numbers in the range
$1,2,\ldots,100$ (see \cite{Borwein:Ferguson:2003}).

We have already mentioned in Section \ref{prelim:sec} that the
periodic Golay number $v>1$ must be even and a sum of two
squares. The following two important necessary conditions for Golay and periodic Golay numbers have been proved more than twenty years ago.

\begin{theorem} \label{Eliahou-Kervaire-Saffari-thm}
{(\rm Eliahou-Kervaire-Saffari \cite{EKS:JCT-A:1990})}
A Golay number is not divisible by any prime
$p\equiv 3 \pmod{4}$.
\end{theorem}

\begin{theorem} \label{Arasu-Xiang-thm}
{(\rm Arasu-Xiang \cite[Corollary 3.6]{Arasu:Xiang:1992})}
If $v=p^t u > 1$ is a periodic Golay number, $p\equiv 3 \pmod{4}$
is a prime number, and $(p,u) = 1$ then $u \geq 2p^{t/2}$.
\end{theorem}
(Since $v$ is a sum of two squares, the exponent $t$ is an
even integer.)

By applying this theorem to the integers in the range
$1,2,\ldots,500$, we deduce that the numbers
$18,36,98,162,242,324,392,484,490$ are not periodic Golay
although each of them is even and a sum of two squares.

It has been shown very recently \cite[Theorem 1]{GSDK:Bio:2014} that the product of a Golay number and a periodic Golay number is again a periodic Golay number. More precisely, the authors of 
that paper show how to ``multiply'' an ordinary Golay pair of length $g$ and a periodic Golay pair of length $d$ to obtain a periodic Golay pair of length $gd$. Consequently, there are infinitely many periodic Golay numbers which are not of the form 
$2^\alpha10^\beta26^\gamma$. We point out that this ``multiplication theorem'' is an easy consequence of \cite[Theorems 13 and 14]{KS:1999}.

As noted in the abstract, there are currently only seven periodic Golay numbers for which we know that they are not Golay numbers. More specifically we have that:
\begin{itemize}

\item
The first periodic Golay pair whose length, $34$, is not a
Golay number has been found in 1998 by {\Dbar}okovi{\'c}
\cite{Djokovic:DesCodes:1998}.
(In fact two non-equivalent such pairs were found.)

\item
Periodic Golay pairs of length $50$ have been found by
{\Dbar}okovi{\'c} \cite{Djokovic:AnnComb:2011} and Kotsireas and
Koukouvinos \cite{KK:2008}.

\item
Periodic Golay pairs of length $58$ have been found by
{\Dbar}okovi{\'c} and Kotsireas \cite{DK:compression:2013}.

\item
Several periodic Golay pairs of length $68$ have been constructed recently in \cite{DKRS::2014}. Such pairs can be constructed also 
by using the method described in \cite{GSDK:Bio:2014}.

\item
Periodic Golay pairs of length $72$ have been constructed recently \cite{DK::2014}.

\item
For examples of periodic Golay pairs of length $74$ see the next section. 

\item
Periodic Golay pair of length $82$ has been found in 2008 by Vollrath \cite{Vollrath:2008} (see also the next section).

\end{itemize}

\section{New periodic Golay pairs}

We have constructed several new periodic Golay pairs, and we
deduce that $74,122,164,202,226$ are periodic Golay numbers. As $74<100$ we know that $74$ is not a Golay number. It is apparently unknown whether $122,164,202,226$ are Golay numbers.
We also construct periodic Golay pairs of length $82$, not
equivalent to the known one \cite{Vollrath:2008}.
In all cases below, the non-equivalence is established by the
method described in \cite{Djokovic:AnnComb:2011}.

The SDSs $(X,Y)$ listed below are given by using the following
compact notation. The parameter set is $(v;r,s;\lambda)$ with
$n=r+s-\lambda$ and $v=2n$. In each case, we make use of a nontrivial subgroup $H$ of the group of units $\bZ_v^\star$ of the ring $\bZ_v$. This subgroup acts on $\bZ_v$ by multiplication
modulo $v$. The orbit of $H$ containing $j\in\bZ_v$ is given by
$H\cdot j=\{hj \pmod{v}:h\in H\}$. The base blocks $X$ and $Y$ are composed of orbits of $H$. Thus we have
$$
X=\bigcup_{j\in J} H\cdot j, \quad Y=\bigcup_{k\in K} H\cdot k,
$$
where $J,K\subseteq\bZ_v$. Instead of listing the elements of
$X$ and $Y$ we shall list only their index sets $J$ and $K$,
respectively. We note that, as a representative $j$ of an orbit
$H\cdot j$ we always choose the smallest integer in that orbit.

A short description of the method we used to construct the
required SDSs and some computational details are given in
section \ref{Algorithm}.

\subsection{Periodic Golay pairs of length $74$}

The parameter set is $(74;36,31;30)$, with $n=37$.
Using the subgroup $H=\{1,47,63\}$ of $Z_{74}^\star$, we give two non-equivalent solutions:

\JK{1,4,6,7,9,12,22,23,28,29,34,42}{1,2,4,6,9,12,17,21,22,37,55}
\JK{1,2,3,6,7,21,22,23,28,29,34,55}{2,4,5,7,9,10,17,21,34,37,42}

\subsection{Periodic Golay pairs of length $82$}

The parameter set is $(82;45,36;40)$, with $n=41$.
Using the subgroup $H=\{1,37,51,57,59\}$ of $Z_{82}^\star$, we give two non-equivalent solutions:

\JK{1,2,11,12,15,17,22,23,30}{1,4,10,12,17,22,23,41}
\JK{1,2,3,6,8,12,17,23,30}{3,5,6,12,17,22,30,41}

\subsection{Periodic Golay pairs of length $122$}

The parameter set is $(122;56,55;50)$, with $n=61$.
Using the subgroup $H =\{1,9,81,95,119\}$ of $Z_{122}^\star$,
we give one solution:

\JK{1,3,6,8,10,13,16,21,23,25,52,61}
{3,4,6,7,13,19,24,25,46,51,52}

\subsection{Periodic Golay pairs of length $164$}

The parameter set is $(164;81,73;72)$, with $n=82$.
Using the subgroup\\
$H=\{1,37,57,133,141\}$ of $Z_{164}^\star$,
we give three non-equivalent solutions:

\JK{4,5,6,10,11,12,16,20,23,25,30,33,46,51,60,65,123}
{3,4,11,13,16,19,20,23,30,33,41,44,46,53,66,82,123}
\JK{1,3,4,6,10,12,13,16,22,23,25,33,34,39,44,46,123}
{2,4,5,8,10,11,12,16,17,33,34,39,41,51,65,82,123}
\JK{4,5,6,8,11,13,20,22,24,30,33,34,39,43,44,65,123}
{2,8,12,13,20,23,25,30,39,41,43,46,51,60,65,82,123}

Additional periodic Golay pairs of length $164$ can be constructed by ``multiplying'' Golay pairs of length $2$ 
with the known periodic Golay pairs of length $82$.

\subsection{Periodic Golay pair of length $202$}
\label{Slucaj:202}

The parameter set is $(202;100,91;90)$, with $n=101$.
Using the subgroup\\
$H=\{1,87,95,137,185\}$ of $Z_{202}^\star$, we give one
solution:

\JK{2,4,9,11,12,13,18,20,22,24,25,26,38,41,50,51,53,55,67,76}
{1,3,4,6,8,9,11,12,16,17,20,25,39,41,48,52,67,76,101}

\subsection{Periodic Golay pairs of length $226$}

The parameter set is $(226;106,105;98)$, with $n=113$.
Using the subgroup\\
$H=\{1,49,109,129,141,143,219\}$ of $Z_{226}^\star$, we give two
non-equivalent solutions:

\JK{1,3,4,5,6,9,10,15,16,36,40,41,43,78,99,113}
{5,8,12,13,15,21,22,24,26,33,34,40,43,78,99}

\JK{3,7,12,13,16,18,20,21,22,40,41,43,55,78,99,113}
{1,2,3,5,10,13,24,26,36,39,40,41,43,78,99}

\section{Algorithm description}
\label{Algorithm}

The algorithm that we used to find the new periodic Golay pairs
is a straightforward adaptation of the algorithm we used in
\cite{DK:JCD:2012} to construct D-optimal matrices.
First we select a subgroup $H$ of $\bZ_v^\star$ and enumerate the
orbits of its action on $\bZ_v$. To construct an SDS with
parameters $(v;r,s;\lambda)$ we first generate two files of
subsets $X$ and $Y$ of size $r$ and $s$, respectively,
of $\bZ_{v}$ such that the corresponding binary sequences $A$ and $B$ pass the PSD test, i.e. satisfy the inequalities
$$
\psd_A(i)\le2v,\quad \psd_B(i) \leq 2v,\quad i=1,\ldots,v/2.
$$
See \cite{DK:JCD:2012,DK:compression:2013} for the precise definition of the PSD function. The subsets $X$ and $Y$ are
constructed as suitable unions of the orbits of the action
described above.
Subsequently we look for a match in the two files, i.e., for two subsets $X$ and $Y$ of size $r$ and $s$ such that the
the pair $(X,Y)$ satisfies the condition stated in
Definition \ref{def-sds}.

Here are some more specific computational details pertaining to the solution for $v=202$ shown in \ref{Slucaj:202}. First we ran
a program for $7$ days to generate a list of about $30$ million
subsets $X$ of $\bZ_{202}$ of size $100$. All of
these sets were made up of $20$ $H$-orbits each of size $5$.
Another $7$-day run of the same program generated about $26$
million subsets $Y$ of size $91$. Each of these subsets was the
union of $18$ orbits of size $5$ and the singleton orbit
$\{101\}$. In both runs we collected only the subsets for which
the corresponding binary sequences pass the PSD test, see
\cite{DK:JCD:2012,DK:compression:2013}. For each of the sets, say $X$, we recorded in a separate file the multiplicities of the
nonzero differences $a-b \pmod{202}$ with $a,b\in X$.
Since these searches were not exhaustive, we applied the transformations $X\to h\cdot X$, $h\in H$, to the output of the first run and for each of the resulting sets we recorded the difference multiplicities. This resulted in a much bigger file containing about $300$ million cases.
In one of the files we replaced the multiplicities $m$ with
$\lambda-m=90-m$, and then ran a program to find matching lines
in the two multiplicity files. The search produced only
two matches but they gave equivalent SDSs. Thus we obtained only
one solution.

\section{Closing comments}

As ordinary Golay pairs are also periodic Golay pairs,
the Golay numbers are also periodic Golay numbers. The known
Golay numbers are exactly the integers $v$ admitting the
factorization $v=2^\alpha 10^\beta 26^\gamma$,
where the exponents $\alpha,\beta,\gamma$ are nonnegative
integers. However, as mentioned earlier, there exist infinitely many periodic Golay numbers which do not admit such factorization. Since all Golay numbers in the range $1,2,\ldots,100$ are known \cite{Borwein:Ferguson:2003}, we deduce that $34,50,58,68,72,74,82$ are the only periodic Golay numbers for which we are presently sure that they are not Golay numbers.

By using the construction (\ref{Konstrukcija}) and the 
``multiplication theorem'' mentioned in section \ref{Known},
we deduce the existence of Hadamard matrices of order 
$2gv$, where $g$ is a Golay number and $v$ a periodic Golay number.

If $v>1$ is a periodic Golay number then $v$ is even, it is a
sum of two squares and satisfies the Arasu-Xiang condition
of Theorem \ref{Arasu-Xiang-thm}. We list all numbers in the
range $1,2,\ldots,300$ which satisfy these three necessary conditions and for which the question whether they are periodic Golay numbers remains open:
$$
90,106,130,146,170, 178,180,194,212,218, 234,250,274,290,292, 298.
$$
This list may be useful to readers interested in constructing new periodic Golay pairs or finding new periodic Golay numbers.

The results of the preprint of this note \cite{DK:arXiv:2013}
(posted on the arXiv), have been already used in  
\cite{GSDK:Bio:2014} for the construction of orthogonal and nearly orthogonal designs for computer experiments.

\section{Acknowledgements}
The authors wish to acknowledge generous support by NSERC.
This work was made possible by the facilities of the Shared
Hierarchical Academic Research Computing Network (SHARCNET) and
Compute/Calcul Canada.

\end{document}